# COMPRESSIBLE MODULES

**Abhay K. Singh**

**Department of Applied Mathematics, Indian School of Mines**

**Dhanbad-826004**

**India**

**Abstract**

The main purpose of this paper is to study under what condition compressible modules are critically compressible. A sufficient condition for the injective hull of a critically compressible module to be critically compressible is also provided. Furthermore we prove sufficient conditions for a critically compressible module to be continuous. In addition, some characterization of critically compressible modules in terms of CS modules, nonsingular modules and cyclic modules are also provided.

**Introduction**

Throughout all rings have non-zero identity elements and all modules are unital right $R-$ modules. Any terminology not defined hare may be found in McConnell and Robson [19] and Gooderl and field [18]. J. Zelmanowitz considered compressible modules in detail in a series of papers [8, 9, 10, 11 and 12]. Zelmanowitz introduced the concept of a weakly primitive ring as a ring having a faithful compressible module. Following [8] the right R-module M is called compressible if for each non-zero submodule $N$ of $M$, there exists a monomorphism $\alpha: M \to N$. In [13], author consider a problem due to Zelmanowitz and study under what condition a uniform compressible module whose non-zero endomorphisms are monomorphisms is critically compressible and they have provided positive answer to this problem for the class of nonsingular modules, quasi-projective modules. In [13], it was proved that the concepts of compressible and critically compressible modules are equivalent over commutative rings or duo rings.

Zelmanowitz [8] claimed that a "compressible uniform module whose nonzero endomorphisms are monomorphisms would be critically compressible". In [2] it was proved that for a commutative ring compressible and critically compressible modules are equivalent and the author called the above statement the "Zelmanowitz Conjecture". Zelmanowitz also raises a question that under what conditions a compressible uniform module whose non-zero endomorphisms are monomorphisms is critically compressible.

Khuri [14] introduced the notion of retractable modules (slightly compressible modules) and gave some results for non-singular retractable modules when the endomorphism ring is (quasi-) continuous. Khuri considered retractable modules in details in series of papers [14, 15, and 16].

This paper is organized in two sections. In section 1 we give preliminary definitions and some results that help me to characterize critically compressible modules.

**1. Preliminaries**

Throughout this paper, it is assumed that R is an associative ring with an identity element. Unless otherwise indicated modules are unitary left modules and homomorphisms are written as right operators. If N is a submodule of M, we write N $\leq$ M and if N is an essential submodule of M then we write N $\trianglelefteq$ M.

Now we will list some definitions and results about critically compressible modules used in this paper.

**Definition 1.1.** Let M be nonzero R-module. Then:

(1) M is called compressible if it can be embedded in each of its nonzero submodules.

(2) M is called critically compressible if it is compressible and, additionally, it cannot be embedded in any of its proper factor modules.

**Definition1.2.** A partial endomorphism of M is a homomorphism from a submodule of M into M.

**Definition1.3.** An R-module M is called uniform if for every pair of non zero submodule $N_1$ and $N_2$, $N_1 \cap N_2 \neq 0$. Clearly M will be uniform if every non zero submodule of M is essential submodule.

**Example1.1.** If every nonzero partial endomorphism on M is monomorphism then M will be uniform but converse is not true. For exp, $Z(p^\infty)$ is uniform Z-module. If we define $f \in End_Z(Z(p^\infty))$ by $f([x]) = [px]$. Then f is nonzero but not monomorphism so, in $Z(p^\infty)$ every nonzero partial endomorphism is not monomorphism.

**Definition1.4.** Given two R-module M and N, M is called N-injective if for every submodule $N_1$ of N, any homomorphism $\alpha : N_1 \to M$ can be extended to a homomorphism $\beta : N \to M$. M is called quasi injective if M is M – injective.

Consider the following condition for a module M.

($c_1$)   Every submodule of M is essential in a direct summand of M.

($c_2$)   If a submodule of M is isomorphic to a direct summand of M, then it is itself a direct summand of M.

($c_3$)   If A & B are direct summand of M with $A \cap B = 0$ then $A + B$ is also a direct summand of M.

M is called CS module (extending module) if it satisfies ($C_1$). M is called continuous module if it satisfies $C_1$ & $C_2$ and quasi continuous if it satisfies $C_1$ & $C_3$.

**Definition1.5.** A nonzero R-module M is called self similar if every non zero submodule of M is isomorphic to M.

**Definition 1.6.** An R-module M is said to be rational extension of $M_1$, in case for each submodule B, $M \supseteq B \supseteq M_1$, $f \in $ Hom (B, M) satisfies $f(M_1) = 0$ if and only if $f = 0$. The injective hull of M denoted by E (M).

**Proposition 1.1.** [7, Prop.1.1] The following conditions are equivalent for a compressible module M:

(1) M is critically compressible;

(2) Every non-zero partial endomorphism on M is monomorphism.

An R-module M is called polyform if every essential submodule of M is dense in M and M is called if every submodule of M is essential. An R-module M is slightly compressible if $Hom_R(M, X) \neq 0$ for every non-zero submodule X of M. Clearly every compressible module is slightly compressible module

**Proposition 1.2.** If every non-zero endomorphism of a slightly compressible module M is monomorphism, then M is compressible.

Since every endomorphism of M is also a partial endomorphism of M, then Proposition 1.1 can be re-written as;

**Proposition 1.3.** Let M is slightly compressible R-module. The following conditions are equivalent:

(i) M is critically compressible;

(ii) Every non-zero partial endomorphism of M is monomorphism.

**2. Compressible and Critically compressible modules**

In this section, we are giving some useful results concerning critically compressible and compressible modules.

**Proposition 2.1.** Every compressible module over a commutative ring R is torsion free.

**Proof:** Since every compressible module over commutative ring is critically compressible, then every non zero partial endomorphism on M is monomorphism.

Take any nonzero partial endomorphism on M. $f_\alpha : M \to M$. Suppose M is not torsion free. Take a torsion element $x \in M$. $\alpha \neq o \in R$ then $\alpha x = 0$ Now define $f_\alpha(x) = \alpha x$ is homomorphism then $f_\alpha(x) = 0 = \alpha x \Rightarrow x \in \text{Ker } f$, but due to critically Compressibility of M, $f_\alpha$ will be monomorphism, then $x = 0$ which is a contradiction. Therefore M is torsion free.

**Proposition2.2.** Let M be compressible module. If M is cyclic and torsion free then M will be critically compressible.

**Proof:** Let N be any nonzero submodule of M and f be any non zero partial endomorphism from N to M. Then $f(n) \neq 0$ for some $o \neq n \in N$. Since otherwise $f(N) = f(Rn) = 0$, a contradiction. Now let $\alpha n \in \text{Ker} f \Rightarrow f(\alpha n) = 0 \Rightarrow \alpha f(n) = 0$ since M is torsion free and $f(n) \neq 0 \Rightarrow \text{Ker} f = 0$ Hence every nonzero partial endomorphism on M is monomorphism. Then M is critically compressible.

**Proposition2.3.** If R be a torsion free ring with unity then $R^R$ is compressible module if and only if $R^R$ it is critically compressible module.

**Proof:** Let $R^R$ is compressible module. Now $R = R.1$, then R is cyclic. Since R is also Torsion free. Then from Prop. 2.2, $R^R$ is critically compressible, converse is obvious.

**Proposition2.4.** Every critically compressible module M is indecomposable.

**Proof:** Suppose $M = M_1 \oplus M_2$ with $0 \neq M_1, M_2 \neq M$. Let $p_i$ be the projection map from M to $M_i$ for $i = 1, 2$. Due to compressibility of M, $p_i$ is monomorphism then Ker $p_i = m_j = 0$ for $i \neq j$ which is a contradiction Hence M is indecomposable.

**Proposition2.5.** Let M be a compressible module. Then every indecomposable injective module M whose nonzero partial endomorphism have uniform non singular kernels then M will be critically compressible.

**Proof:** Let N be a nonzero submodule of M. Now Let $f : N \to M$ be any non zero partial endomorphism on M. suppose Ker $f \neq 0$ since Ker f is uniform, N = Ker f [by 2, Prop. 2.2] then N is non singular. Therefore N is rational extension of Ker f [by 3,

Prop 5, p. 59]. Hence f (Kerf) = 0 $\Rightarrow$ f = 0 which is a contradiction. Hence Ker f = 0 then every non zero partial endomorphism on M is monomorphism. Then M is critically compressible.

**Proposition2.6.** Every critically compressible module is CS (extending module) and it is continuous if and only if every partial monomorphism is isomorphism.

**Proof:** Let M be a critically compressible module. Then every non zero partial endomorphism on M is monomorphism. So M is uniform. Hence every critically compressible module is CS module (extending module).

Let f : M $\to$ M be any non zero partial endomorphism on M. Due to critically compressibility of M, f will be monomorphism. Suppose f is not surjective then f (M) C M and f (M) will be proper essential submodule of M. Now define g : M $\to$ f (M) by g (x) = f (x), then g is isomorphism. Since M is continuous, f (M) being on isomorphic copy of M cannot be a proper essential submodule of M, a contradiction. Hence f is surjective. Converse is obvious.

We note that, every non-zero submodule of compressible or critically compressible module is also compressible or critically compressible module respectively.

It is clear that if E(M) is compressible or critically compressible then M is also compressible or critically compressible module respectively. Now we are giving sufficient condition for an injective hull of compressible or critically compressible module to be compressible or critically compressible.

**Proposition2.7.** Let M is a critically compressible then E(M) will be critically compressible if

(1) M is stable under partial endomorphism of E (M)

(2) E (M) is rational extension of M.

**Proof:** Let M is critically compressible. Now suppose $M^1$ be any non zero submodule of E (M) and h : $M^1$ $\to$ E (M) be any zero partial endomorphism on E (M). Suppose

Ker h ≠ 0. Then Ker h will be submodule of E (M) ⇒ M ∩ Ker h ≠ 0. Now f = h/N is any non zero partial endomorphism on M. Since M is stable under endomorphism of E (M). Also f ≠ 0 since E (M) is rational extension of M. Thus Ker f = Kerh ∩M ≠ 0 which is a contradiction. Hence E (M) is critically compressible.

Now if M is non-singular then rational extension of M and injective hull of M coincides and this leads immediately to the following.

**Corollary2.8.** Let M and E (M) be compressible module and M be nonsingular module stable under partial endomorphism of E (M). Then E(M) is critically compressible module if and only if M is compressible module.

A ring with unit 1 is said to be an S.P. ring if Z (M) is a direct summand of M for every R-module M, where Z (M) = {x ∈ M | R Δ 0 (x)}, 0 (x) denotes the annihilator ideal of x. *In* [7] author claimed that if M is compressible module then M is singular or nonsingular. Now we are giving here the result that for a S.P. ring every critically compressible module is either singular or nonsingular.

**Proposition2.9.** Let R be a S.P. ring. Then any critically compressible module over R is either singular or non singular.

**Proof.** Let M is neither singular nor non singular, then Z (M) is a non zero direct summand and proper submodule of M. Now by prop. 2.4 M is indecomposable Hence Z (M) = M, a contradiction. Hence M is either singular or non singular.

**Proposition2.10.** A ring R with unit $1 \neq 0$ is self similar and in which every non zero endomorphism is monomorphism if and only if it is a principal ideal ring without zero divisor.

**Proof.** Assume R is self similar and in which every non zero endomorphism is monomorphism. Then R has no zero divisor. To prove R is a principal ideal ring, Let S be any non zero ideal of R. Take $\alpha \in R$, r≠0 then there exist a isomorphism f : R →

Rr and g : Rr → S. Let   f (1) = xr and g (xr) = b, Now if y∈S then y = g of (a)  for some a∈R

y=g o f (a.1)

= g (a f(1))

= a g (xr) = ab ∈R

Hence , S = Rb is a principal ideal.

The converse is obvious.

**Corollary 2.11.** If R is a principal ideal ring without zero divisor then $R^R$ is critically compressible.

**Proof.** If R is Principal ideal ring without zero divisor then R is self similar in which every non zero endomorphism is monomorphism (By prop. 2.10). Then from [1, thereon 4.1], $R^R$ is critically compressible.

Hence the following implications hold :

*R is self similar in which every non zero endomorphism is monomorphism* ⇒ *R is a principal ideal ring without zero divisor* ⇒ $R^R$ *is critically compressible.*

**Proposition 2.12.**  Let M be compressible module and every non zero f ∈ $Hom_R$ (M, M') is monomorphism, where M' is quasi injective hull of M. Then M is critically compressible module.

**Proof.** Let g : N → M be any non zero partial endomorphism on M. Such that $f_{|N} = g$ by assumption f is monomorphism and so g is monomorphism. Hence every non zero partial endomorphism on M is monomorphism. Then M is critically compressible.

**Proposition 2.13.** Let R be a noetherian ring and R be a module over ring R. Then following conditions are equivalent:

    (i)      $R^R$ is compressible module

    (ii)     $R^R$ is critically compressible

    (iii)    $Monomorphism_R$ (R, U) is non zero for every uniform submodule U of R.

(iv) Monomorphism$_R$ (R, U) is non zero for every cyclic uniform submodule U of R.

**Proof.** (1) $\Rightarrow$ (2) from [1, theorem 3.2]

(2) $\Rightarrow$ (3) $\Rightarrow$ (4) is clear.

Now (4) $\Rightarrow$ (1)

Let N be any non zero submodule of R. Let m be any non zero element of N. By hypothesis mR is noetherian and hence mR contains a uniform submodule U of R. Let u be any non zero element of U. Then uR is cyclic uniform submodule of R. By (4) mono$_R$ (R, UR) is non zero. Hence mono$_R$ (R, N) is non zero. Hence $R^R$ is compressible module.

**References:-**